\documentclass[12pt]{article}

\usepackage{amsfonts}

\usepackage{amssymb}

\usepackage{amsmath}

\usepackage{eufrak}

\usepackage[dvips]{graphicx}

\usepackage{latexsym}

\usepackage{verbatim}

\newcommand{\ind}{\makebox[1em]{\raisebox{-.5ex}[0ex][0ex]{\makebox[0em]%

{$\smile$}}\raisebox{.4ex}[0ex][0ex]{\makebox[-.02em]{$|$}}}}

\newcommand{\dep}{\makebox[1em]{\raisebox{.3ex}[0ex][0ex]%

{$\not$}\makebox[.7em]{\ind}}}

\newcommand{\mdep}{\makebox[1em]{\raisebox{1.5ex}[0ex][0ex]{\makebox[0em]%

{$\scriptscriptstyle m$}}\makebox[.6em]{$\dep$}}}

\newcommand{\mind}{\makebox[1em]{\raisebox{1.5ex}[0ex][0ex]{\makebox[0em]%

{$\scriptscriptstyle m$}}\makebox[.6em]{$\ind$}}}

\newcommand{\bi}{\begin{itemize}}

\newcommand{\ei}{\end{itemize}}

\DeclareMathOperator{\im}{Im}

\DeclareMathOperator{\ann}{Ann}

\newtheorem{theorem}{Theorem}[section]

\newtheorem{lemma}[theorem]{Lemma}

\newtheorem{fact}[theorem]{Fact}

\newtheorem{claim}{Claim}

\newtheorem{claim1}{Claim}

\newtheorem{subclaim}{Subclaim}

\newtheorem{corollary}[theorem]{Corollary}

\newtheorem{proposition}[theorem]{Proposition}

\newtheorem{definition}[theorem]{Definition}

\newtheorem{remark}[theorem]{Remark}

\newtheorem{conjecture}[theorem]{Conjecture}

\newtheorem{question}[theorem]{Question}

\newtheorem{example}[theorem]{Example}

\title{Locally finite profinite rings}

\author{Jan Dobrowolski\footnote{Research supported by NCN Grant no. 2012/05/N/ST1/02850} \hspace{1mm} and Krzysztof Krupi\'nski}

\date{}

\begin{document}

\maketitle

\begin{abstract} We investigate the structure of locally finite profinite rings. We classify (Jacobson-) semisimple locally finite profinite rings as products of complete matrix rings of bounded cardinality over finite fields, and we prove that the Jacobson radical of any locally finite profinite ring is nil of finite nilexponent. 
Our results apply to the context of small compact $G$-rings, where we also obtain a description of possible actions of $G$ on the underlying ring.

\end{abstract}



\section{Introduction}

By a locally finite ring we mean a ring whose every finitely generated subring is finite and by a profinite ring we mean the limit of any inverse system of finite rings. 

Now, we recall the notions of small compact $G$-groups and $G$-rings, which form important subclasses of the class of small Polish structures defined in \cite{1}.


\begin{definition}
Let $G$ be a Polish group. A compact $G$-group [$G$-ring] is a pair $(R,G)$, where $R$ is a compact topological group [ring] and $G$ acts on $R$ continuously as a group of automorphisms. We say that $(R,G)$ is small, if for every $n<\omega$ there are only countably many orbits under the action of $G$ on $R^{\times n}$ (the set of $n$-tuples).
\end{definition}

Our general motivation is to understand the structure of groups and rings in the context of small Polish structures. A particularly interesting and accessible situation is the case of small compact $G$-groups and $G$-rings. The initial motivation for the current work was to describe the structure of small compact $G$-rings. Under the additional assumption of $nm$-stability, we 
know a lot by \cite{2,7}, e.g. \cite[Theorem 3.2]{7} tells us that small, $nm$-stable compact $G$-rings are nilpotent-by-finite, and it is 
conjectured that they are null-by-finite which was confirmed if the so-called ${\cal NM}$-rank of the ring is less than $\omega$. However, without $nm$-stability, not much is known. We know for example that that the Cartesian power $R^\omega$ of any finite ring $R$, considered together with the Polish group $G$ of all permutations of coordinates, is a small compact $G$-ring which, of course, does not need to contain nilpotent elements.

By Fact \ref{fact1}, we know that each small compact $G$-ring is a locally finite profinite ring. So, our goal is to describe the structure of locally finite profinite rings admitting a structure of a small compact $G$-ring. It turns out, however, that our proofs work in the general context of locally finite profinite rings, or even more generally, in the context of profinite rings whose all 1-generated subrings are finite. So, the main structural results of this paper are stated in this general context, without using the notion of small compact $G$-rings, except Corollary \ref{action of G}, where we additionally describe the action of $G$ on the ring $R$.

\begin{definition}
We say that a ring $R$ is weakly locally finite, if every 1-generated subring of $R$ is finite.
\end{definition}

In order to describe the structure of a ring $R$, it is important to understand the structure of the semisimple ring $R/J(R)$ and of the radical ring $J(R)$, where $J(R)$ is the Jacobson radical of $R$. The main results of this paper are:  a complete classification of semisimple, weakly locally finite profinite rings established in Theorem \ref{thm1}, and important information on the structure of the Jacobson radical of weakly locally finite profinite rings obtained in Theorem \ref{thm2}.




\begin{theorem}\label{thm1}
Let $R$ be a topological ring. Then, $R$ is a semisimple, weakly locally finite profinite ring if and only if $R$ is isomorphic (as a topological ring) to a direct product of complete matrix rings over finite fields with only finitely many non-isomorphic rings occurring as factors in this product.
\end{theorem}

From this classification, one immediately concludes that semisimple, weakly locally finite profinite rings coincide with semisimple, locally finite profinite rings.

The above theorem not only yields a complete classification of the class of semisimple locally finite profinite rings, but also of the class of semisimple rings admitting a structure of a small compact $G$-ring (see Corollary \ref{class});
in Corollary \ref{action of G}, we also describe possible actions of $G$ on the ring in question.

A problem which is ``complementary'' to the description of semisimple rings from a given class is the problem of describing the Jacobson radical of rings belonging to that class. Our second main result is the following (see also Corollary \ref{loc fin implies nil}).

\begin{theorem}\label{thm2}
If $R$ is a weakly locally finite profinite ring, then $J(R)$ is nil of finite nilexponent. More generally, each nil profinite ring has finite nilexponent.
\end{theorem}

In particular, the Jacobson radical of each small compact $G$-ring is nil of finite nilexponent.
This generalizes a similar result proved in \cite{5} for so called small profinite rings. Recall that
a small profinite ring [group] is the limit $R$ of a countable inverse system of finite rings [groups] together with a closed subgroup $Aut^*(R)$ of the group of automorphisms respecting the distinguished inverse system such that $Aut^*(R)$ has only countably many orbits on
$n$-tuples for all $n<\omega$. In particular, every small profinite ring [group] is a small compact $G$-ring [$G$-group] (with $G:=Aut^*(R)$).

Recall that it was proved \cite{5} that the Jacobson radical of a small profinite ring is open, and the following conjectures were formulated.

\begin{conjecture}\label{nilpotent}
The Jacobson radical of a small profinite ring $R$ is nilpotent. In particular, $R$ has an open nilpotent ideal.
\end{conjecture}

\begin{conjecture}\label{null}
A small profinite ring has an open null ideal.
\end{conjecture}

These conjectures are interesting in their own rights, but additional motivation standing behind them comes from \cite{5} and \cite{7}. To explain this, recall the main conjecture concerning small profinite groups.

\begin{conjecture}\label{small profinite groups are abelian}
A small profinite group has an open abelian subgroup.
\end{conjecture}

By \cite[Corollary 2.4]{7} (see also \cite[Corollary 3.13]{7}), we know that Conjecture \ref{small profinite groups are abelian} implies Conjecture \ref{null} which, of course, implies Conjecture \ref{nilpotent}. On the other hand, \cite[Theorem 3.5]{5} tells us that Conjecture \ref{nilpotent} for commutative rings implies an important intermediate conjecture towards the proof of Conjecture \ref{small profinite groups are abelian}, namely that each small soluble profinite group has an open nilpotent subgroup; using \cite[Theorem 2.10]{7}, we get that Conjecture \ref{null} implies that each small soluble profinite group has an open abelian subgroup.

Several reductions of Conjecture \ref{nilpotent} for commutative rings were obtained in \cite{5}. In Section \ref{remarks}, we prove some further reductions of that conjecture. The main one (see Proposition \ref{polynomials}) roughly says that if Conjecture \ref{nilpotent} is false, then there is a counter-example which generically does not satisfy any polynomial identities which are not satisfied by obvious reasons.

At the end notice that the counterpart of the second part of Conjecture \ref{nilpotent} for small compact $G$-rings is false by Theorem \ref{thm1}. As to the first part, we do not know.

\begin{question}\label{q}
Is it true that the Jacobson radical of each small compact $G$-ring is nilpotent?
\end{question}

It is clear, however, that the Jacobson radical of a small compact $G$-ring does not need to be null-by-finite. For this, take the Cartesian power $R^\omega$ of any non-trivial finite nil ring $R$ and consider it with the Polish group $G$ of all permutations of coordinates. This is a nilpotent, small compact $G$-ring which is not null-by-finite.

\section{Preliminaries}\label{prel}

In this paper, we always assume rings to be associative, but we do not assume them to be commutative or unital. By an ideal we mean a two-sided one. By a topological ring we mean a ring equipped with a Hausdorff topology under which multiplication, addition and additive inversion are continuous functions.

Recall some basic notions from ring theory. An element $r$ of a ring $R$ is nilpotent of nilexponent $n$ if $r^n=0$ and $n$ is the smallest number with this
property. The ring is nil [of nilexponent $n$] if every element is nilpotent [of nilexponent $\leq n$ and there is an element of nilexponent
 $n$]. The ring is nilpotent of class $n$ if $r_1\cdots r_n=0$ for all $r_1,\ldots,r_n\in R$ and $n$ is the smallest number 
with this property. An element $r$ is null if $rR=Rr=\{0\}$. The ring is null if all its elements are. 
If $r \in R$, then $\ann(r) :=\{ a \in R: ar=ra=0 \}$ is the two-sided annihilator of $r$ in $R$. Note that $\ann(r)$ is always a subgroup of $R^+$, and if $R$ is commutative, then $\ann(r)$ is an ideal.

Now, we recall fundamental issues concerning Jacobson radicals; for more details see \cite{14}.
The Jacobson radical of a ring $R$, denoted by $J(R)$, is the collection of all
elements of $R$ satisfying the formula $\phi(x)=\forall y\exists z(yx+z+zyx=0)$ (that is, it is the set of all elements which generate quasi-regular left ideals, i.e. left ideals consisting of left quasi-regular elements which are defined as those elements $x \in R$ for which there is $z \in R$ such that $x+z+zx=0$). Equivalently, $J(R)$ is the unique largest quasi-regular left [or right] ideal. Another equivalent definition says that $J(R)$ is the intersection of all the maximal regular left [or right] ideals, where a left ideal $I$ is said to be regular if there is $a \in R$ such that $x-xa \in I$ for all $x \in R$ (notice that for rings with $1$ all ideals are regular). 
For any ring $R$, $J(R)$ is a two-sided ideal. We say that $R$ is semisimple if $J(R)=\{0\}$. $R/J(R)$ is always a semisimple ring.  Every nil left [or right] ideal is contained in $J(R)$; in particular, if $R$ is nil, then $J(R)=R$.

For any left ideal $I$ of $R$, define $m_I$ to be the largest (two-sided)  ideal of $R$ contained in $I$. Then, for any maximal left regular ideal $I$, $R/m_I$ is a (left) primitive ring (i.e. a ring having a left faithful irreducible module, namely the module $R/I$). The next fact follows from \cite[Theorem 2.1.4]{14}. 

\begin{fact}\label{matrix}
If $I$ is a maximal left regular ideal of a ring $R$ and the quotient $R/m_I$ is finite, then $R/m_I$
 is the complete matrix ring $M_k(F)$ over a finite field $F$.
\end{fact}


\begin{remark}\label{za+a}
Let $R$ be any ring. For every $x \in J(R)$ and $a \in R$ such that $xa + a=0$, we have $a=0$.
\end{remark}

{\em Proof.} 
Since $x \in J(R)$, we get that there is some $z$ such that $zx+z+x=0$. Then $0=zxa+za+xa=-za+za-a=-a$. 
\hfill $\square$

\begin{corollary}\label{loc fin implies nil}
If $R$ is a weakly locally finite ring, then $J(R)$ is nil.
\end{corollary}

{\em Proof.} Take any $x \in J(R)$. By assumption, there are $n>m\geq 1$ such that $x^n=x^m$. Then $(-x^{n-m})x^m + x^m=0$ and $-x^{n-m} \in J(R)$. Thus, $x^m=0$ by Remark \ref{za+a}. \hfill $\square$\\

The following remark follows easily from the definition of $J(R)$.

\begin{remark}\label{J(R) closed}
If $R$ is a compact topological ring, then $J(R)$ is a closed ideal.
\end{remark}

The following fact \cite[Proposition 5.1.2]{riza} yields a
characterization of when a topological unital ring is profinite.

\begin{fact}\label{classical} 
Let $R$ be a topological unital ring (so, in particular, a Hausdorff
topological space). Then the following conditions are
equivalent:\begin{enumerate}
\item $R$ is profinite, i.e.\ the inverse limit of finite rings.
\item $R$ is compact.
\item $R$ is compact and totally disconnected.
\item $R$ is compact and there is a basis of open neighbourhoods of $0$
consisting of open ideals.\end{enumerate}
\end{fact}

Now, we will give some basic information about small compact $G$-rings.
Let $(R,G)$ be a compact $G$-ring. For any finite $C\subseteq R$, by $G_C$ we denote the pointwise stabilizer of $C$ in $G$, and for a finite tuple $a$ of elements of $R$, by $o(a/C)$ we denote the orbit of $a$ under the action of $G_C$ (and we call it the orbit of $a$ over $C$).
Then we have that $(R,G)$ is small iff for every finite $C\subseteq R$ there are only countably many orbits on $R$ over $C$.

\begin{fact}\label{fact1}
Every small compact $G$-ring $R$ is locally finite and profinite.
\end{fact}
{\em Proof.} Local finiteness of $R$ follows as in \cite[Proposition 5.7]{1}. Namely, consider any finite subset $S$ of $R$, and let $\overline{\langle S \rangle}$  be the closure of the subring generated by $S$. As each element of $\overline{\langle S \rangle}$ is fixed by $G_S$, smallness implies that $\overline{\langle S \rangle}$ is countable. But it is also a compact group, so it must be finite.

By Fact \ref{classical}, we see that $R$ is profinite when $R$ has a unit, and we will use this information to show that $R$ is always profinite. By the Baire category theorem and local finiteness of $R$, we get that for some non-zero $n<\omega$ the set $\{r\in R : nr=0\}$ has non-empty 
interior. Hence, we can cover $R$ with finitely many tranlates of this set, which yields that $R$ has a finite characteristic $c$.
Put $R_1=R\times \mathbb{Z}_c$, and define $+$ and $\cdot$ 
on $R_1$ by $(a,k)+(b,l)=(a+b,k+_cl)$ and $(a,k)\cdot(b,l)=(ab+l\times a+k\times b,k\cdot_cl)$, 
where $+_c$ and $\cdot_c$ are addition and multiplication modulo $c$, and $l\times a:=a+\dots+a$ ($l$-many times).
Then, $R_1$ is an unital compact ring, and we can treat $R$ as a clopen ideal of $R$ in the natural way. By Fact \ref{classical}, $R_1$ is profinite, so also $R$ is profinite.
\hfill $\square$\\

The following remarks provide some examples of small compact $G$-rings, and they will be used later.
\begin{remark}\label{power}
Let $R$ be a finite ring. 
Consider the action of the permutation group $G=S_{\omega}$ on $R^{\omega}$ given by 
$(\sigma\cdot f)(i)=f(\sigma^{-1} (i))$. Then, $(R^{\omega},G)$ is a small compact $G$-ring.
\end{remark}


\begin{remark}\label{product}
Let $(R_1,G_1)$, $(R_2,G_2)$ be two small compact $G$-rings. Then also $(R_1\times R_2,G_1\times G_2)$ is a small compact $G$-ring, where the action is given by $(g_1,g_2)\cdot(r_1,r_2)=(g_1\cdot r_1,g_2\cdot r_2)$.
\end{remark}

Suppose we have a profinite ring $R$ which is the inverse limit of a
distinguished countable inverse system, or equivalently, $R$ is a compact
topological ring with a distiguished countable basis of open neighbourhoods of $0$
consisting of clopen ideals.

\begin{definition} A profinite ring regarded as profinite 
structure is a pair of the form $(R,Aut^*(R))$, where $Aut^*(R)$ is a
closed subgroup of the group $Aut^0(R)$ of all automorphisms of $R$
respecting the inverse system defining $R$ (equivalently, 
$Aut^*(R)$ is a closed subgroup of the group of all automorphisms of
$R$ fixing setwise the clopen ideals from the distiguished basis of
open neighbourhoods of 0).
The group $Aut^*(R)$ is called the structural group of $R$ and
$Aut^0(R)$ is the standard structural group of $R$.
\end{definition}

\begin{definition} A profinite ring $(R,Aut^*(R))$  is small if there are only countably many
orbits under $Aut^*(R)$ on finite tuples over $\emptyset$ (equivalently, if it is small regarded as a compact $G$-ring). 
\end{definition}

For the general definition of a profinite structure see
\cite{4}. 
Until the end of this section, by a profinite ring we mean a profinite ring regarded a profinite structure. 
Let $(R,Aut^*(R))$ be a profinite ring. 
As usual, we will write $R$ having in mind either the topological ring $R$ or the pair $(R,Aut^*(R))$.

 By an $A$-invariant subset of $R^{\times n}$ we mean a subset invariant under the pointwise stabilizer of $A$ (denoted by
$Aut^*(R/A)$), and by an $A$-closed subset we mean a closed
and $A$-invariant subset of $R^{\times n}$; if we do not want to specify $A$,
we write $*$-invariant or  $*$-closed. For $a \in R^{\times n}$ and
$A\subseteq R$ we define $o(a/A)= \{f(a):f\in Aut^*(R/A) \}$, the 
orbit of $a$ over $A$.

The following remark follows directly from the definition.

\begin{remark}\label{1.one} A profinite ring $R$ has a descending chain $(I_n:n<\omega)$
of open $\emptyset$-invariant ideals forming a basis of open
neigbourhoods of $0$ (hence with trivial intersection)
\end{remark}

It is clear than whenever $I$ is an $A$-closed ideal of $R$, then $R/I$ can be treated as a profinite ring (regarded as profinite structure) with the structural group induced by $Aut^*(R/A)$.

Since every orbit is a closed subset of $R$ (as a continuous image of a compact space $Aut^*(R/A))$, it follows, by the Baire category theorem, that if $R$ is small, then over any finite subset $A$ there exists an open orbit.

By $acl(A)$ we denote the set of elements of $R$ which have
finite orbits over $A$. 
We say that $g\in R$ is generic over $A$ if $o(g/A)$ is open; $g$ being generic means that it is generic over $\emptyset$ or over a set of parameters which is obvious from the context. We
say that $a$ is $m$-independent from $b$ over $A$ (denoted
$a \mind_A b$) if the orbit $o(a/Ab)$ is open in $o(a/A)$. Otherwise
$a$ is $m$-dependent on $b$ over $A$ (written $a \mdep_A b$). The
relation $\mind$ enjoys some properties analogous to the ones that forking
independence has in simple theories.

The following fact is due to Newelski in \cite{3}.

\begin{fact}\label{new} Let $(G,Aut^*(G))$ be a small profinite group, and $A$ a
finite subset.
An $A$-invariant subgroup of $G$ is $A$-closed. The group generated
by any family of $A$-invariant sets is $A$-closed, and generated in
finitely many steps from finitely many sets. There is no infinite
increasing chain of $A$-invariant subgroups of $G$. In particular, all
characteristic subgroups of $G$ are $\emptyset$-closed.  
\end{fact}

It was observed in \cite{5} that in order to obtain Conjecture \ref{nilpotent} for commutative rings, it is enough to prove it only for rings with the characteristic and nilexponent equal to the same prime number (this relies on the Nagata-Higman theorem which yields that nil rings of nilexponent smaller than the characteristic are nilpotent). Also, the following reduction was obtained there (see \cite[Proposition 4.4]{5}).

\begin{fact}\label{reduction} To show Conjecture \ref{nilpotent} for commutative rings, one can assume that for each non-zero $a\in R$, the quotient ring $R/\ann(a)$ is not nilpotent. 
\end{fact}

\section{Main Results}

This section is devoted mainly to the proofs of our main results about weakly locally finite profinite rings. 
We start from a certain characterization of Jacobson radicals in profinite rings.

\begin{proposition}\label{additional}
Let $R$ be a profinite ring, and let $\cal I$ be the family of all maximal regular left open (so clopen) ideals of $R$. Recall that for any left ideal $I$, $m_I$ denotes the unique largest (two-sided) ideal of $R$ contained in $I$.
\begin{enumerate}
\item Each $I \in \cal I$ is a maximal regular left ideal of $R$.
\item For each $I \in \cal I$, $m_I$ is the unique largest (two-sided) clopen ideal of $R$ contained in $I$; in particular, $R/m_I$ is finite.
\item $J(R)= \bigcap_{I \in \cal I} m_I=\bigcap \cal I$.
\end{enumerate}
\end{proposition}
%
{\em Proof.}
(1) Take $I \in {\cal I}$. Then $I$ is of finite index in $R$. So, any left ideal extending $I$ is also clopen.\\
(2) Take any $I \in \cal I$. Since it is clopen, there is a clopen (two-sided) ideal $m$ of $R$ contained in $I$. Then $m_I +m$ is a clopen ideal of $R$ contained in $I$, so, by the maximality of $m_I$, we get that $m_I=m_I+m$ is clopen.\\
(3) Put $J_1(R)= \bigcap \cal I$ and $J_2(R)= \bigcap_{I \in \cal I} m_I$. By \cite[Theorem 1.2.1]{14}, we know that $J(R) = \bigcap_{J \in \cal J} m_J$, where $\cal J$ is the collection of all maximal regular left ideals. So, by (1), we get get $J(R) \subseteq J_2(R) \subseteq J_1(R)$. It remains to show that $J_1(R) \subseteq J(R)$. By the proof of \cite[Theorem 1.2.2]{14}, we see that it is enough to show that for any $x \in J_1(R)$ the set $A:=\{yx+y : y \in R\}$ is all of $R$. Suppose for a contradiction that $A$ is proper. Since $R$ is compact, $A$ is a regular left closed ideal of $R$. Thus, there exists a clopen ideal $I$ such that $A+I$ is a proper regular left clopen ideal of $R$. Hence, $A$ is contained in a maximal regular left clopen ideal $J$. Then, for every $y \in R$, we have $x \in J_1(R) \subseteq J$ and $yx+y \in J$, so $y \in J$. Thus, $J=R$, a contradiction.
\hfill $\square$\\

In the proof of Theorem \ref{thm1}, we will use this proposition only for weakly locally finite unitial rings. In this case, one can give a direct proof (not referring to the proof of \cite[Theorem 1.2.2]{14}), which we do below for the reader's convenience.\\

\noindent 
{\em Proof of Proposition \ref{additional}(3) for weakly locally finite unitial rings.}
As before, it is enough to show that $J_1(R) \subseteq J(R)$, and for this, it is enough to show that $J_1(R)$ is nil. Since $R$ is weakly locally finite, we will be done if we show that $1-x$ is left invertible in $R$ for all $x\in J_1(R)$. So, take any $x\in J_1(R)$ and suppose for a contradiction that $1\notin R(1-x)$.
By compactness, $R(1-x)$ is closed in $R$, so there is an open ideal $I$ of $R$ such that $1\notin R(1-x)+I$.
We can extend $R(1-x)+I$ to a maximal left open ideal $J$ of $R$. Then $x\in J_1(R)\subseteq J$ and $1-x\in J$, so $1\in R$. Thus, $J=R$, a contradiction.
\hfill $\square$\\

Recall that for a filed $F$, $M_k(F)$ denotes the ring of all matrices of size $k \times k$ with entries from $F$.
\begin{lemma}\label{bound}
Let R be a weakly locally finite profinite ring. Then there exists $n<\omega$ such that for every finite field $F$ and a number $k<\omega$ if there exists an epimorphism $f:R\to M_k(F)$, then $k<n$ and $|F|<n$.
\end{lemma}
{\em Proof.}
For any $1<m<\omega$, consider the polynomial $$w_m=\prod_{0<i<j\leq m}(x^i-x^j).$$ Since $R$ is weakly locally finite, we have that $R=\bigcup_{m<\omega}Z_R(w_m)$, where $Z_R(p)$ denotes the set of all zeros  in $R$ of 
a polynomial $p$. Since each $Z_R(w_m)$ is closed in $R$, we get, by the Baire category theorem, that for some $m<\omega$, $Z_R(w_m)$ has a non-empty interior. Take $a\in R$ and an open ideal $I$ of $R$ of index $s$ such that $a+I\subseteq Z_R(w_m)$. 
Since $f[I]$ is an ideal of $M_k(F)$, we get that it is either trivial or equal to $M_k(F)$.



Now, $[M_k(F):f[I]]\leq s$, and if $f[I]=\{ 0\}$, then $[M_k(F):f[I]]\geq 2^{k^2}$. So, either $2^{k^2}\leq s$, or $f[I]=M_k(F)$ and then $f[a+I]=f(a)+f[I]=M_k(F)$. In the latter case, since $(\forall x\in a+I)(w_m(x)=0)$, we get that $(\forall x\in M_k(F))(w_m(x)=0)$. Using the matrix with $1$'s right above the diagonal and $0$'s elsewhere, one easily gets that $k\leq \sum_{i=1}^{m-1} i(m-i)=\frac{m(m-1)(m+1)}{6}$.
So, we obtain the following bound on $k$: 
$$k\leq \max\left(\sqrt{log_2(s)},\frac{m(m-1)(m+1)}{6}\right).$$

Now, if $f[I]=M_k(F)$, then $(\forall x\in F)(w_m(x)=0)$ (as $F$ embeds into $M_k(F)$), so $|F|\leq deg(w_m)$. On the other hand, if
$f[I]=\{0\}$, then $|M_k(F)|\leq s$, so $|F|\leq s$. Thus, in any case, $$|F|\leq \max(s, deg(w_m)).$$
\hfill $\square$\\

Now, we proceed to the proof of Theorem \ref{thm1}.\\

\noindent
{\em Proof of Theorem \ref{thm1}.} The implication $(\leftarrow)$ is easy. Let $R$ be topologically isomorphic to the product of complete matrix rings over finite fields with only finitely many factors up to isomorphism. Clearly $R$ is profinite and locally finite. To see semisimplicity, is is enough to use a classical fact that the complete matrix ring over any field is semisimple (see e.g. \cite[Theorem 1.2.6]{14}) and to show that products of semisimple rings are semisimple, which is a very easy exercise.

Now, we turn to the proof of $(\rightarrow)$.  Let $R$ be our semisimple, weakly locally finite profinite ring.
\begin{claim}\label{cl0}
We can assume that R is unital.
\end{claim}
{\em Proof of Claim \ref{cl0}.}
Suppose that the theorem is true for unital rings.
Define $R_1$ as in the proof of Fact \ref{fact1}.
Since $R$ is semisimple, we get that $R\cap J(R_1)=J(R)=\{0\}$. Thus, $R$ is isomorphic (via the quotient map) to a closed ideal of $R_1/J(R_1)$.
Since $R_1/J(R_1)$ is semisimple, we get, by our assumption, that it is of the form $\prod_{i=1}^m M_{n_i}(F_i)^{\kappa_i}$
for some cardinal numbers $\kappa_i$. It is easy to see that every closed ideal of this ring consists precisely of elements with zeros
on a fixed set of coordinates, so it is also isomorphic to a ring of this form. In particular, $R$ can be presented in this form.
\hfill $\square$\\

Now, assume that $R$ is unital. 
Let $\cal I$ be the family of all maximal left open ideals of $R$. Let ${\cal I}_1$ be the subfamily of those $I \in \cal I$ for which $m_I$ is minimal in the family of all $m_I$'s, $I \in \cal I$. Finally, for each $m \in \{ m_I : I \in {\cal I}_1\}$ choose exactly one ideal $I$ from the set $\{ I \in {\cal I}_1 : m_I =m\}$, and denote the collection of all ideals obtained in this way by ${\cal I}_2$.

By Proposition \ref{additional}(2), for each $I \in \cal I$, $R/m_I$ is finite, so Proposition \ref{additional}(1) and Fact \ref{matrix} imply that $R/m_I$ is isomorphic to a complete matrix ring over a finite field. Thus, using Lemma \ref{bound}, we get that the family $\{ m_I : I \in \cal I\}$ is well-founded with respect to inclusion. Therefore,
$\bigcap_{I \in \cal I} m_I = \bigcap_{I \in {\cal I}_2} m_I$. 
By Proposition \ref{additional}(3) and the assumption that $R$ is semisimple, we end up with
\begin{equation}\tag{$*$}
\bigcap_{I \in {\cal I}_2} m_I=J(R)=\{0\}.
\end{equation}

Define a homomorphism $f:R\to \prod_{I\in {\cal I}_2} R/m_I$ as the diagonal of the quotient homomorphisms. By $(*)$, we get that $f$ is injective. Also, $f$ is continuous, so $\im(f)$ is closed.

Now, consider any $I,J\in {\cal I}_2$ with $I\neq J$. Since $m_J$ is not contained in $m_I$, we get that $(m_I+m_J)/m_I$ is a non-trivial ideal of $R/m_I$. Hence, as we know that $R/m_I$ is a complete matrix ring over a finite field, we get that $(m_I+m_J)/m_I$ is equal to $R/m_I$, so $m_I+m_J=R$. Therefore, since $R$ is unitial, it follows from the Chinese Remainder Theorem that $\im(f)$ is dense in $\prod_{I\in A} R/m_I$. So,  $\im(f)=\prod_{I\in A} R/m_I$, and $f$ is an isomorphism. Moreover, by Lemma \ref{bound}, the complete matrix rings $R/m_I$ have bounded size. This completes the proof of Theorem \ref{thm1}.
\hfill $\square$

\begin{corollary}\label{class}
Let $\cal{A}$ denote the class all topological rings isomorphic to a product of complete matrix rings of bounded size over finite fields, and let ${\cal A}_0$ be its subclass consisting of products of only countably many matrix rings.
\begin{enumerate}
\item $\cal A$ is the class of all semisimple weakly locally finite profinite rings.
\item $\cal A$ is the class of all semisimple locally finite profinite rings.
\item ${\cal A}_0$ is the class of all semisimple topological rings admitting a structure of a small compact $G$-ring.
\end{enumerate}
\end{corollary}
{\em Proof.} (1) is a restatement of Theorem \ref{thm1}.\\ 
(2) follows from (1) and the observation that each ring from $\cal A$ is locally finite.\\
(3) By Remarks \ref{power} and \ref{product}, any member of ${\cal A}_0$ admits a structure of a small compact $G$-ring. By Fact \ref{fact1} and (2), we get that each semisimple topological ring $R$ admitting a structure of a small compact $G$-ring belongs to $\cal A$. However, \cite[Proposition 3.9]{1} tells us that $R$ is second countable. Thus, we conclude that $R \in {\cal A}_0$. 
\hfill $\square$\\

Using this classification, we can also describe possible actions of $G$ on $R$ for semisimple small compact $G$-rings $(R,G)$.
\begin{corollary}\label{action of G}
Let ${\cal A}$ be defined as above. Take any $R\in \cal{A}$. Present $R$ in the form $\prod_{i=1}^m M_{n_i}(F_i)^{\kappa_i}$ 
for some cardinal numbers $\kappa_i$ so that the rings $M_{n_i}(F_i)$ are non-isomorphic for distinct $i$'s. Then the group of all topological automorphisms of $R$ is equal to $$\{\prod_{i=1}^m f_{\sigma_i}\circ g: \sigma_i\in Sym(\kappa_i),g\in \prod_i Aut(M_{n_i}(F_i))^{\kappa_i}\},$$ where $f_\sigma (x)(\alpha)=x(\sigma^{-1}(\alpha))$. Hence, it is isomorphic to the semidirect product of $\prod_i Sym(\kappa_i)$ and $\prod_i Aut(M_{n_i}(F_i))^{\kappa_i}$.  Thus, if $(R,G)$ is a small compact $G$-ring, then we can treat $G$ as a subgroup of the above group, acting on $R$ in the natural way.
\end{corollary}
{\em Proof.} 
Consider an arbitrary topological automorphism $f$ of $R$. Fix any $i$ and consider any $\alpha \in \kappa_i$.

We have that $f[\{r\in R:(\forall (j,\beta)\neq (i,\alpha))( r(j,\beta)=0)\}]$ is a closed ideal of $R$ isomorphic to $M_{n_i}(F_i)$, and since it is not a product of two non-trivial rings, we get that it is equal to $\{r\in R:(\forall (j,\beta)\neq (i,\gamma)) (r(j,\beta)=0)\}$ for some $\gamma\in \kappa_i$. Define $$\sigma_i(\alpha)=\gamma.$$

By composing $f$ with canonical isomorphisms $M_{n_i}(F_i)\rightarrow\{r\in R:(\forall (j,\beta)\neq (i,\alpha))( r(j,\beta)=0\})$ and $\{r\in R:(\forall (j,\beta)\neq (i,\gamma))( r(j,\beta)=0)\}\rightarrow M_{n_i}(F_i)$, we obtain an automorhism $g_{i,\alpha}$ of $M_{n_i}(F_i)$. Put 
$$h=\prod_{i=1}^m (f_{\sigma_i}\circ \prod_{\alpha\in\kappa_i}g_{i,\alpha}).$$

By the choice of $h$, we get that $h$ agrees with $f$ on elements of $R$ having one-element supports. Since $f$ and $g$ are algebraic isomorphisms, we get that they agree on all elements with finite supports. Finally, we conclude, by continuity of $f$ and $h$, that $f=h$. This completes the proof.
\hfill $\square$

\begin{question}
Let $(R,G)$ be a semisimple small compact $G$-ring. By Corollary \ref{class}, we 
know that $R \in {\cal A}_0$, and we have that the group $G$, treated as a permutation group of $R$, is a subgroup of the concrete Polish group $Aut(R)$ described in Corollary \ref{action of G}. Is it the case that $G$ is also a topological subgroup of $Aut(R)$, i.e. is the topology on $G$ inherited from $Aut(R)$? Equivalently, is $G$ a closed subgroup of $Aut(R)$? 
\end{question}

Now, we turn to the proof Theorem \ref{thm2}.\\

\noindent
{\em Proof of Theorem \ref{thm2}.} If $R$ is a weakly locally finite profinite ring, then $J(R)$ is nil by Corollary \ref{loc fin implies nil} and closed by Remark \ref{J(R) closed}, so $J(R)$ is a nil profinite ring. Thus, the second part of the theorem is indeed more general than the first part. So, from now on, assume that $R$ is a nil profinite ring. Then $J(R)=R$.

Since $R$ is nil, we get, by the Baire category theorem, that there is a non-zero $n<\omega$ such that the set $\{r \in R:r^n=0\}$ has a non-empty interior.
Take an open ideal $I$ of $R$ and a coset $e=a+I$ such that $e\subseteq  \{r \in R:r^n=0\}$.

We put 
$Z_{-1}=\{ 0 \}$,
and define inductively $$Z_{i+1}=\{r\in R:IrI/(Z_i\cap IrI) \;\, \mbox{is nilpotent}\}.$$ 
It is easy to see that $Z_i$ is an ideal of $R$ for each $i$ (if $r_1 \in Z_{i+1}$ yields the quotient $Ir_1I/(Z_i\cap IrI)$ of nilpotency class $c_1$ and $r_2$ yields the quotient of nilpotency class $c_2$, then $r_1+r_2$ yieds the quotient of nilpotency class at most $c_1+c_2$).

\begin{claim1}\label{cl1}
For all $k\in \{0,\dots, n-1\}$, $y \in I$ and $x\in e$ we have $x^{n-k}yx^{n-k}\in Z_{k-1}$.
\end{claim1}
{\em Proof of Claim \ref{cl1}.}
We proceed by induction on $k$. The case $k=0$ is trivial. Suppose the claim is true for a certain $k<n-1$. 
Then $w^{n-k}\in Z_k$ for every $w\in e$. So, in $R/Z_k$, we have
$$\begin{array}{ll}
0+Z_k=(x+yx^{n-k-1})^{n-k}+Z_k=\\
=x^{n-k-1}yx^{n-k-1}+x^{n-k-2}(yx^{n-k-1})^2+\dots +(yx^{n-k-1})^{n-k}+Z_k=\\
=x^{n-k-1}yx^{n-k-1}+zx^{n-k-1}yx^{n-k-1}+Z_k,
\end{array}$$ 
where $z=x^{n-k-2}y+x^{n-k-3}yx^{n-k-1}y+\dots +(yx^{n-k-1})^{n-k-2}y$.
But $z+Z_{k}\in R/Z_{k}= J(R/Z_{k})$, because $R/Z_{k}$ is nil.
Hence, by Remark \ref{za+a}, $x^{n-k-1}yx^{n-k-1}+Z_k=0+Z_k$, which shows that the claim is true for $k+1$.
\hfill $\square$\\

Applying Claim \ref{cl1} for $k=n-1$, we easily get that $e\subseteq Z_{n-1}$, so $Z_{n-1}$ is clopen.
For $x\in R$ and $\overline{i}=(i_1,\dots, i_{2m})$ (where each $i_j$ is from $R$), we define 
$$\overline{i}*x = i_1xi_2i_3xi_4\dots i_{2m-1}xi_{2m}.$$
\begin{claim1}\label{cl2}
The following statement is true for all $k \in \omega$.
$$\begin{array}{cc}
(\forall x\in Z_{k-1})(\exists m_1<\omega)( \forall \overline{i_1}\in I^{\times 2m_1})(\exists m_2<\omega)(\forall \overline{i_2}\in I^{\times 2m_2})\dots\\
 (\exists m_k<\omega)
( \forall \overline{i_k}\in I^{\times 2m_k})
(\overline{i_k}*(\overline{i_{k-1}}*(\dots (\overline{i_1}*x)\dots ))=0).
\end{array}$$
\end{claim1}
{\em Proof of Claim \ref{cl2}.} 
The claim follows from the definition of the $Z_{j}$'s by induction on $k$. The case $k=0$ is trivial. For the induction step, 
consider any $x \in Z_{k-1}$. Let $m_1$ be the nilpotency class of $IxI/(Z_{k-2} \cap IxI)$. Then, for any $\overline{i_1} \in I^{\times 2m_1}$ we have that $\overline{i_1} * x \in Z_{k-2}$, so the assertion follows from the inductive hypothesis. \hfill $\square$\\

In the next claim, we will show that in the statement from Claim \ref{cl2} one can move all the existential quantifiers to the left, obtaining a statement of the form $(\forall x \in Z_{k-1})\exists \forall$.

\begin{claim1}\label{cl3}
The following statement is true for all $k \in \omega$.
$$\begin{array}{cc}
(\forall x\in Z_{k-1})(\exists m_1,\dots,m_k <\omega)( \forall \overline{i_1}\in I^{\times 2m_1},\dots,\overline{i_k}\in I^{\times 2m_k})\\
(\overline{i_k}*(\overline{i_{k-1}}*(\dots (\overline{i_1}*x)\dots ))=0).
\end{array}$$
\end{claim1}
{\em Proof od Claim \ref{cl3}.} By induction on $k$, we will show that whenever $x \in R$ is such that 
$$\begin{array}{cc}
(\exists m_1<\omega)( \forall \overline{i_1}\in I^{\times 2m_1})(\exists m_2<\omega)(\forall \overline{i_2}\in I^{\times 2m_2})\dots
 (\exists m_k<\omega)
( \forall \overline{i_k}\in I^{\times 2m_k})\\
(\overline{i_k}*(\overline{i_{k-1}}*(\dots (\overline{i_1}*x)\dots ))=0),
\end{array}$$
then 
$$\begin{array}{cc}
(\exists m_1,\dots,m_k <\omega)( \forall \overline{i_1}\in I^{\times 2m_1},\dots,\overline{i_k}\in I^{\times 2m_k})\\
(\overline{i_k}*(\overline{i_{k-1}}*(\dots (\overline{i_1}*x)\dots ))=0).
\end{array}$$
%
This together with Claim \ref{cl2} will finish the proof.

The cases $k=0$ and $k=1$ are trivial, as there are no quantifiers to switch. Now, take $k>1$ and assume that the statement is true for numbers less than $k$. Consider any $x \in R$ satisfying the assumption of our statement. By the inductive hypothesis, we get 
\begin{equation}\tag{$*$}\label{eq *}
\begin{array}{cc}
(\exists m_1<\omega)( \forall \overline{i_1}\in I^{\times 2m_1})(\exists m_2,\dots,m_k <\omega)( \forall \overline{i_2}\in I^{\times 2m_2},\dots,\overline{i_k}\in I^{\times 2m_k})\\
(\overline{i_k}*(\overline{i_{k-1}}*(\dots (\overline{i_1}*x)\dots ))=0).
\end{array}
\end{equation}

Now, the goal is to switch $( \forall \overline{i_1}\in I^{\times 2m_1})$ with $(\exists m_2,\dots,m_k <\omega)$. We will do this in $2m_1$ steps, switching at every step all existential quantifiers $\exists m_2, \dots, \exists m_k$ with one universal quantifier corresponding to one of the variables in the sequence $\overline{i_1}$. We will only show how to switch all these existential quantifiers with the universal quantifier corresponding to the last variable in $\overline{i_1}$, as the other steps can be done in a similar fashion.

Denote $\overline{i_1}=(s_1,\dots,s_{2m_1})$, and fix $s_1,\dots,s_{2m_1-1} \in I$. Put $$t=s_1xs_2s_3xs_4\dots s_{2m_1-1}x.$$ 
For any $\overline{m}=(m_2,\dots,m_k) \in (\omega \setminus \{ 0 \})^{k-1}$ define
$$D_{\overline{m}}= \{ i \in I : ( \forall \overline{i_2}\in I^{\times 2m_2},\dots,\overline{i_k}\in I^{\times 2m_k})
(\overline{i_k}*(\overline{i_{k-1}}*(\dots *(ti)\dots ))=0)\}.$$

By (\ref{eq *}), we have $I = \bigcup_{\overline{m}} D_{\overline{m}}$. It also follows that each $D_{\overline{m}}$ is a closed subset of $I$. Hence, by the Baire category theorem, there is some $\overline{m_0}$ such that $D_{\overline{m_0}}$ has a non-empty interior in $I$. Thus, for some 
$a_1,\dots,a_w\in I$ we have that 
$$I=(a_1+D_{\overline{m_0}})\cup\ldots\cup (a_w+D_{\overline{m_0}}).$$ 
We also know that there are $\overline{m_1},\dots,\overline{m_w}$ such that $a_1\in D_{\overline{m_1}},\dots, a_w\in D_{\overline{m_w}}$. So, in order to finish the proof, it is enough to show the following subclaim (in which $\overline{m_i}$'s and $a_i$'s are NOT the particular tuples or elements chosen above).

\begin{subclaim}
For any $l \geq 1$ and   $\overline{m}=(m_1,\dots,m_l) \in (\omega \setminus \{ 0 \})^{l}$ define
$$D_{\overline{m}}' = \{ a \in I : ( \forall \overline{i_1}\in I^{\times 2m_1},\dots,\overline{i_l}\in I^{\times 2m_l})
(\overline{i_l}*(\overline{i_{l-1}}*(\dots  *(\overline{i_1}*a)\dots ))=0)\}.$$ 
Let $r \geq 1$. Then, for every $\overline{m_1},\dots,\overline{m_r}$ there is $\overline{m}$ such that for any
 $a_1 \in D_{\overline{m_1}}', \dots, a_r \in D_{\overline{m_r}}'$ one has $a_1 + \ldots +a_r \in D_{\overline{m}}'$.
\end{subclaim}

%
%
\noindent
{\em Proof of Subclaim 1.} The proof is by induction on $l$.
Consider the base step $l=1$. Take any $\overline{m_1},\dots,\overline{m_r}<\omega$.
Let $m=\overline{m_1}+\dots +\overline{m_r}$. Consider any $a_1 \in D_{\overline{m_1}}', \dots, a_r \in D_{\overline{m_r}}'$. Then, for every $\overline{i_1} \in  I^{\times 2m}$ the element $\overline{i_1}*(a_1+\ldots+a_r)$ is a sum of elements from the sets $(Ia_1I)^{\overline{m_1}}, \dots, (Ia_rI)^{\overline{m_r}}$ which are all equal to $\{ 0 \}$, so  $\overline{i_1}*(a_1+\ldots+a_r)=0$ and $a_1+\ldots+a_r \in D_{\overline{m}}'$.

The induction step is similar. Take any $\overline{m_1}=(m_1^1,\dots,m_l^{1}), \dots,  \overline{m_r}=(m_1^r,\dots,m_l^{r})$. Let $m'=m_1^1+\dots +m_1^r$.
Consider any $a_i \in D_{\overline{m_i}}'$ for $i=1,\dots,r$. Then, the elements from the set 
$$\{ \overline{i_1}*(a_1+\ldots+a_r) : \overline{i_1} \in I^{\times 2m'}\}$$ 
are sums of a bounded number of elements from the sets $$\{ \overline{i_1}*a_1 : \overline{i_1} \in I^{\times 2m_1^1}\}, \dots, \{ \overline{i_1}*a_r : \overline{i_1} \in I^{\times 2m_1^r}\},$$ and we finish using the inductive hypothesis.\hfill $\square$\\

So, we have proved the statement formulated at the beginning of the proof of Claim \ref{cl3} which together with Claim \ref{cl2} completes the proof of Claim \ref{cl3}.\hfill $\square$\\


Consider the statement from the last claim for $k=n$. Using the fact that $Z_{n-1}$ is compact, we can apply the same trick as in the proof of Claim \ref{cl3} to switch all the existential quantifiers with the quantifier $\forall x \in Z_{n-1}$, and so we get that there are $m_1,\dots,m_n<\omega$ such that  
$$(\forall x\in Z_{n-1})( \forall \overline{i_1}\in I^{\times 2m_1},\dots,\overline{i_k}\in I^{\times 2m_n})\\
(\overline{i_n}*(\overline{i_{n-1}}*(\dots (\overline{i_1}*x)\dots ))=0).$$
Hence, $I\cap Z_{n-1}$ has finite nilexponent. Since $R/(I\cap Z_{n-1})$ is finite, it has also finite nilexponent. 
Thus, we conclude that $R$ has finite nilexponent, which completes the proof of Theorem \ref{thm2}.
\hfill $\square$\\

Summarizing, we have the following corollary of Theorems \ref{thm1} and \ref{thm2}.
\begin{corollary}
Every weakly locally finite profinite ring is (nil of finite nilexponent)-by-(product of complete matrix rings over finite fields with only finitely many factors up to isomorphism).
\end{corollary}

Having Theorem \ref{thm2}, one could ask if we can strengthen it by showing that the Jacobson radical of a locally finite profinite ring is nilpotent. The following easy example shows that this is not always true.

\begin{example}\label{series}
Let $p$ be a prime number and let $F_0$ be the free commutative nil ring of nilexponent $p$ and of
characteristic $p$ on generators $(x_i:i<\omega)$. For $n<
\omega$ let $I_n$ be the ideal generated by $\{x_i:i\ge n\}$. Then each
quotient ring $F_0/I_n$ is finite and nilpotent; their inverse limit, say $F$, is the
free commutative profinite ring which is nil of nilexponent $p$ and
of characteristic $p$ with free topological generators $\{y_n:n<\omega\}$, where
$y_n=(x_n+I_k : k\in \omega)$. We see that $F=J(F)$ is locally finite, but it is not nilpotent.
\end{example}

However, for the class of small compact $G$-rings that question (see Question \ref{q}) is open. In particular, we do not know whether the ring $F$ from Example \ref{series} considered together with $G$ being the group of all topological automorphisms of $F$ is small. We should remark here that by an easy counting argument (see \cite[Example 1 in Section 5]{5}), one can check that $F$ cannot be a counter-example to conjecture \ref{nilpotent}, i.e., it does not admit a structure of a small profinite ring.

















We finish with some observations concerning small profinite rings.

\begin{proposition}\label{products}
\begin{enumerate}
\item Suppose $(\prod_{i\in I} R_i, G)$ is a small profinite ring. Then only finitely many $R_i$'s are not null rings.
\item Suppose $(\prod_{i\in I} H_i, G)$ is a small profinite group. Then only finitely many $G_i$'s are non-abelian 
groups.
\end{enumerate}
\end{proposition}
{\em Proof.}
(1)
Put $J=\{i\in I:R_i$ is not null$\}$. For any $i\in J$ choose $s_i\in R_i$ such that the two-sided annihilator of $s_i$ in $R_i$ is not equal to $R_i$.
Let $o$ be an open orbit in $\prod_{i\in I} R_i$. Then, for some finite $I_0\subseteq I$ and elements $r_i\in R_i$, $i\in I_0$, we have that $\{r\in\prod_{i\in I} R_i:(\forall i\in I_0)( r(i)=r_i)\}\subseteq o$. Define $x,y\in\prod_{i\in I} R_i$ by: $x(i)=r_i$ if $i\in I_0$ and $x(i)=0$ 
otherwise; $y(i)=r_i$ if $i\in I_0$, $y(i)=0$ if $i\in I\backslash (I_0\cup J)$ and $y(i)=s_i$ if $i\in J\backslash I_0$. Since $x,y\in o$ and the two-sided annihilator of $x$ is open in $R$, we get that the same is true about $y$. This clearly implies that $J$ is finite.\\
(2) is similar, using centralizers instead of annihilators.
\hfill $\square$\\

Part (2) of the above proposition slightly strengthens Remark 4.3 from \cite{3}, where it is additionally assumed that the structural group consists of automorphisms respecting the inverse system $H_1\leftarrow H_1\times H_2\leftarrow\dots$. Part (1) is a ring counterpart of (2).

Notice that Theorem 3.1 from \cite{5} follows easily from Theorems \ref{thm1}, \ref{thm2} and Proposition \ref{products}(1).

\section{Remarks on Conjecture \ref{nilpotent}}\label{remarks}

In this section, we prove some reductions for Conjecture \ref{nilpotent}. It has already been recalled in Section 1 that in order to prove Conjecture \ref{nilpotent} for commutative rings, one can assume that the ring in question is nil of nilexponent $p$ and of characteristic $p$ for some prime number $p$, which justifies this assumption in the results below.

Throughout this section, we will skip the structural group $Aut^*(R)$.
\begin{lemma}\label{powers}
Fix any prime number $p$. Let $R$ be a non-nilpotent, commutative small profinite ring of characteristic $p$ which is nil of nilexponent $p$.
Then, for every $g\in G$ generic over $\emptyset$ and for every $j\in \{1,2,\dots,p-1\}$, we have $g^j\neq 0$.
\end{lemma}
{\em Proof.}
We proceed by induction on $j$. The conclusion is clear for $j=1$, so suppose $j>1$ and that the lemma holds for smaller numbers.
Suppose for a contradition that $g^j=0$ for some generic $g\in R$. Then there is an open ideal $I$ of $R$ such that $(g+i)^j=0$ for all $i\in I$. 
Therefore, $i^jg^{j-1}=g^{j-1}(g+i)^j=0$ for $i\in I$. Take any $i\in I$ which is $m$-independent from $g$. If $R/\ann(i^j)$ is not nilpotent,
then it is a commutative small profinite ring of charatceristic $p$ and nilexponent $p$ (by the Nagata-Higman theorem),  and since $g+\ann(i^j)$ is its generic satisfying $(g+\ann(i^j))^{j-1}=0+\ann(i^j)$, we would get a contradition with the inductive hypothesis. So, $R/\ann(i^j)$ is nilpotent. Thus, taking $I':=\{i\in I: i\mind g\}$ and $K:=\{b\in R: R/\ann(b)$ is nilpotent$\}$, we get that $\{i^j: i\in I'\}\subseteq K$. By the claim in the proof of Proposition 4.4 in \cite{5}, we have that $K$ is a nilpotent ideal of $R$.
It follows from Fact \ref{new} that $K$ is also closed. Since $I'$ is a dense subset of $I$, we get, by the continuity of the mapping $x\mapsto x^j$, that $\{i^j:i \in I\}\subseteq K$. Thus, $I/(I\cap K)$ is a nil ring of nilexponent not greater than $j<p$, whence we get, by the Nagata-Higman theorem, that it is nilpotent. Since also $K$ is nilpotent, we get that $I$ is nilpotent. Hence, $R$ is nilpotent, a contradiction.
\hfill $\square$\\

Now, we will make the main observation of this section. We will assume that $R/\ann(a)$ is non-nilpotent for every $a\in R\backslash \{0\}$. This assumption is justified by Fact \ref{reduction}.

\begin{proposition}\label{polynomials}
Let $R$ be as in the lemma, and assume additionally that $R/\ann(a)$ is non-nilpotent for every $a\in R\backslash \{0\}$. 
Then, for every polynomial $f(x_1,\dots, x_n)\in F_p[x_1,\dots,x_n]\backslash \{0\}$ with $deg_{x_1}(f),\dots,deg_{x_n}(f)<p$ and for all independent tuples $(g_1,\dots,g_n)$ of generics of $R$ we have that $f(g_1,\dots, g_n)\neq 0$.
\end{proposition}
{\em Proof.}
We proceed by induction on $n$. 
Suppose that $n\geq1$ and that the proposition is true for smaller positive natural numbers (the argument will also cover the base induction step). We can present any $f$ satysfying the assumptions in the form
$$f=h_0(x_1,\dots,x_{n-1})+h_1(x_1,\dots,x_{n-1})x_n+\dots +h_{p-1}(x_1,\dots,x_{n-1})x_n^{p-1}$$
(where $h_0,\dots,h_{p-1}$ are constants in the case when $n=1$).
We will prove the conclusion by induction on the number $k$ of non-zero polynomials among $h_0,\dots,h_{p-1}$.

Suppose first that $k=1$ and $f(x_1,\dots,x_n)=h_i(x_1,\dots,x_{n-1})x_n^i$ for some $i>0$ (we can assume that $deg_{x_n}f>0$ by the inductive hypothesis). Take any tuple $(g_1,\dots, g_n)$ of independent generics of $R$. Put $a=h_i(g_1,\dots, g_{n-1})$. By the inductive hypothesis of the first induction, $a\neq 0$. By assumptions on $R$, we have that $R/\ann(a)$ is non-nilpotent, so it has characteristic and nilexponent equal to $p$. The coset $g_n+\ann(a)$ is a generic element of this ring, so, by Lemma \ref{powers}, $g_n^i\notin \ann(a)$. Hence, $f(g_1,\dots,g_n)=ag_n^i\neq 0$.

Now, we turn to the inductive step, where we assume that $k>1$. Then, $$f(x_1,\dots,x_n)=h_{i_1}(x_1,\dots,x_{n-1})x_n^{i_1}+h_{i_2}(x_1,\dots,x_{n-1})x_n^{i_2}+\dots +h_{i_k}(x_1,\dots,x_{n-1})x_n^{i_k},$$
where $i_1<i_2<\dots<i_k$ are all indices $i$ for which $h_i$ is non-zero (if $i_1=0$, then by $h_{i_1}(x_1,\dots,x_{n-1})x_n^{i_1}$ we mean just $h_{i_1}(x_1,\dots,x_{n-1})$).
By the inductive hypothesis of the second induction, we get that the element
$$
\begin{array}{ll}
g_n^{p-i_k}f(g_1,\dots,g_n) = h_{i_1}(g_1,\dots,g_{n-1})g_n^{i_1+p-i_k}+h_{i_2}(g_1,\dots,g_{n-1})g_n^{i_2+p-i_k}+\dots +\\
+h_{i_{k-1}}(g_1,\dots,g_{n-1})g_n^{i_{k-1}+p-i_k}+h_{i_k}(g_1,\dots,g_{n-1})g_n^{i_k+p-i_k}=h_{i_1}(g_1,\dots,g_{n-1})g_n^{i_1+p-i_k}+\\ 
+h_{i_2}(g_1,\dots,g_{n-1})g_n^{i_2+p-i_k}+\dots +h_{i_{k-1}}(g_1,\dots,g_{n-1})g_n^{i_{k-1}+p-i_k}
\end{array} 
$$
is non-zero, so $f(g_1,\dots,g_n)\neq 0$.
\hfill $\square$

\begin{corollary}
If Conjecture \ref{nilpotent} is not true in the class of commutative rings, then there is a counter-example, say $R$, to it such that the topological ring $F$ defined in Example \ref{series} (for some prime $p$) is topologically isomorphic to a closed subring of $R$.
\end{corollary}
{\em Proof.}
As it has already been explained, by results of \cite[Section 4]{5}, we can take $R$ to be a counter-example to Conjecture \ref{nilpotent} which satisfies the assumptions (and hence, the conclusion) of Proposition \ref{polynomials}. Let $(I_n)_{n<\omega}$ be a decresing chain of open ideals of $R$ with trivial intersection.
We choose inductively a sequence of independent generics $g_0,g_1,\dots$ of $R$ such that $g_{n+1}\in I_{k_n}$, where $k_n$ is a natural number such that for any polynomial $f(x_0,\dots, x_n)\in F_p[x_0,\dots,x_n]\backslash \{0\}$ with $deg_{x_0}(f),\dots,deg_{x_n}(f)<p$, we have that $f(g_0,\dots,g_n)\notin I_{k_n}$. Let $S$ be the closure of the subring of $R$ generated by $\{g_i:i<\omega\}$.

We define $\phi :F\to S$ as follows. Consider any $y\in F$. Choose polynomials $p_i(t_0,\dots,t_i)$, $i<\omega$, such that 
$y=\lim_i p_i(y_0,\dots,y_i)$ for $y_i$'s defined in Example \ref{series}. Then, by the choice of the sequence $(g_i)$, the sequence $p_i(g_0,\dots,g_i)$ is convergent, and we define $\phi(y)$ to be its limit. It is easy to check that $\phi: F \to S$ is a well-defined isomorphism of topological rings.
\hfill $\square$\\

Note that Conjecture \ref{nilpotent} would be proved if in the above proof we were able to embed topologically $F$ into $R$ as an invariant (under all topological automorphisms, or only under the ones coming from the structural group of $R$) ring. Indeed, if there was such an embedding, then $F$ with the structural group induced by the structural group of $R$ would be a small profinite ring, which is impossible by the comment in the paragraph below Example \ref{series}.

\noindent
{\bf Address:}\\
Instytut Matematyczny, Uniwersytet Wroc\l awski,\\
pl. Grunwaldzki 2/4, 50-384 Wroc\l aw, Poland.\\[3mm]
{\bf E-mail addresses:}\\
Jan Dobrowolski: dobrowol@math.uni.wroc.pl \\
Krzysztof Krupi\'nski: kkrup@math.uni.wroc.pl

\end{document}